\documentclass[11pt,reqno]{article}

\usepackage[LCY,T1]{fontenc}

\usepackage[style=trad-abbrv,isbn=false,url=false,doi=false, firstinits=true]{biblatex}

\usepackage[russian,main=english]{babel}
    \addto{\captionsenglish}{}
\usepackage{csquotes}
\usepackage{amsfonts}
\usepackage{amscd}
\usepackage{color}
\usepackage{todonotes}
\usepackage{amsmath}
\usepackage{float}
\usepackage{amsthm}
\usepackage{amssymb}
\usepackage{mathrsfs}
\usepackage{amstext}
\usepackage{graphicx}
\usepackage{tikz}
\usetikzlibrary{matrix}
\usepackage{bbm}

\usepackage{mathtools}

\numberwithin{equation}{section}

\theoremstyle{plain}
\newtheorem{satz}{Theorem}[section]
\newtheorem{thm}[satz]{Theorem}
\newtheorem{defi}[satz]{Definition}
\newtheorem{cor}[satz]{Corollary}
\newtheorem{lem}[satz]{Lemma}
\newtheorem{prop}[satz]{Proposition}
\newtheorem{rem}[satz]{Remark}

\newcommand{\supp}{\operatorname{supp}}

\newcommand{\mix}{{\rm mix}}

\newcommand{\re}{\ensuremath{\mathbb{R}}}
\newcommand{\R}{\ensuremath{{\re}^d}}

\newcommand{\N}{\ensuremath{\mathbb{N}}}
\newcommand{\n}{\ensuremath{{\N}_0}}

\newcommand{\zz}{\ensuremath{\mathbb{Z}}}
\newcommand{\Z}{{\ensuremath{\zz}^d}}

\newcommand{\C}{\ensuremath{\mathbb{C}}}

\newcommand{\tor}{\ensuremath{\mathbb{T}}}

\newcommand{\id}{\mbox{id}}

\newcommand{\ca}{\ensuremath{\mathcal A}}

\newcommand{\cF}{\ensuremath{\mathcal F}}

\newcommand{\bk}{\ensuremath{\mathbf k}}

\newcommand{\bj}{\ensuremath{\mathbf j}}

\newcommand{\bx}{\ensuremath{\mathbf x}}

\newcommand{\by}{\ensuremath{\mathbf y}}

\newcommand{\hyp}{\operatorname{%
\mathchoice{%
\includegraphics{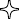}}{%
\includegraphics{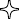}}{%
\includegraphics{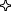}}{%
\includegraphics{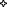}}
}}

\newcommand{\rank}{{\rm rank \, }}

\newlength{\fixboxwidth}
\makeatletter
\def\blfootnote{\xdef\@thefnmark{}\@footnotetext}
\makeatother

\usepackage{hyperref}
\begin{document}
\title{Gelfand and sampling numbers of weighted Wiener spaces}

\author{Moritz Moeller\footnote{
moritz.moeller@math.tu-chemnitz.de} \\\\ Chemnitz University of Technology, Faculty of Mathematics}

\date{November 9, 2026}

\maketitle

\begin{abstract}
The aim of this paper is to establish sharp asymptotic bounds on the Gelfand numbers of weighted Wiener spaces \(\ca^r_\theta\) measured in \(L_q\). This is done by using results for Gelfand numbers of sequence spaces and applying them to each layer of a hyperbolic cross to obtain results for Gelfand numbers between different (weighted) Wiener spaces. These can then be used to establish bounds in \(L_q\) by combining them with known results and techniques for Besov spaces. This significantly narrows a gap in the current research, since the Gelfand numbers provide a lower bound on the (non-linear) sampling numbers, whereas upper bounds on this quantity have been studied in multiple papers via an upper bound given by the best \(m\)-term approximation.
\end{abstract}

\section{Introduction}
The non-linear sampling numbers \(\varrho_m(X)_Y\), originally conceived in \cite{WaWo01}, describe the behavior of the worst case error over all functions from the function class of \(X\) of the reconstruction using standard information (\(m\) point evaluations) with an arbitrary algorithm. They are notoriously difficult to deal with directly and are instead often bounded from below and above by other error terms.

In \cite{JUV23} the authors developed a method that enables one to bound the non-linear sampling numbers in \(L_2\) from above by a sparse approximation from a dictionary plus a residual term, both in \(L_\infty\). In the setting of the trigonometric system these are the best \(m\)-term approximation and the error of best trigonometric approximation.
This result was then generalized to (non-linear) sampling numbers in general \(L_q\) spaces in \cite{MPU24}. Furthermore, the target norm for the best \(m\)-term approximation was replaced by a Wiener norm and an additional factor in \cite{Kri23} and \cite{MSU26}. In many cases this yields results that are at least as good if not better than measuring in \(L_\infty\).
The most recent refinements of this inequality were made in \cite{ABG26}, where the error of best approximation was replaced by some tail bounds. This has positive ramifications for the computational costs for implementations of the underlying methods.

All of this refined machinery has been used to show upper bounds on the sampling numbers of certain function spaces measured in \(L_q\). In particular, results for weighted Wiener spaces \(\ca^r_\theta\), which are the spaces that we also want to consider in this paper, have been derived this way.

These are function spaces defined via a weighted sum of their Fourier coefficients. Their norm is defined as 
 \begin{equation}
     \|f\|_{\ca^{r}_{\theta}} \coloneqq \Big( 
\sum_{\bk\in\Z} \prod_{i=1}^d (1+|k_i|)^{ r\theta}|\hat{f}(\bk)|^\theta \Big)^\frac{1}{\theta}. \end{equation}

Approximation in these spaces has recently been studied in numerous papers. In \cite{NNS22} the authors studied multiple different \(s\)-numbers of mixed Wiener spaces. The best \(m\)-term trigonometric approximation for generalized Wiener spaces was considered in \cite{DPW24}. Further results for generalized Wiener spaces can be found in \cite{Shid26}.  In \cite{KLT23} sparse grid methods in Wiener spaces were studied. The connection to Barron classes was first formally observed in \cite{Vo22} and results for linear sampling numbers were shown in \cite{KPUU23}. There are many more results in these spaces, but core to our topic of interest are \cite{Moe23, MSU26, MSU25}. They deal with the best \(m\)-term approximation for these spaces in \(L_q\) and derive upper bounds for the sampling numbers via the above-mentioned method.
In \cite[Corollary 6.2]{MSU26} the following upper bound is shown this way,
\begin{equation} \label{eq:intro_samp_upper}
 \varrho_{m}\big(\ca_\theta^r \big)_{L_q} \lesssim m^{-( r+{\frac{1}{\theta}+\frac{1}{q}-1})} \log(m)^{(d-1) r + 3 ( r + {\frac{1}{\theta}+\frac{1}{q}-1})}.
\end{equation}

However, there are no corresponding general lower bounds. Indeed the only (almost) sharp lower bound given is from \cite[Theorem 5.1]{Moe23} and states
\begin{equation}
     \varrho_{m}\big(\ca_\theta^r \big)_{L_2} \gtrsim m^{-(r - \frac{1}{2} + \frac{1}{\theta})} \log(m)^{(d-1)r}.
\end{equation}
This is derived from a sharp two-sided bound on the Gelfand numbers. The techniques used in \cite{Moe23} heavily rely on the fact that \(L_2 = \ca_2\) and can therefore transition to the sequence space setting where results on the Gelfand numbers are readily available.

In this paper we want to show lower bounds on the sampling numbers of \(\ca^r_\theta\) in general \(L_q\) with \(1 \leq q \leq \infty\). We will do so by developing sharp results for the Gelfand numbers of \(\id: \ca^r_\theta \to L_q\) which in turn will rely on bounds for the Gelfand numbers of \(\id: \ca^r_\theta \to \ca_\eta \). Finally, these results will be derived from known bounds on the Gelfand numbers between sequence spaces, from, e.g., \cite{FPRU10} or \cite{Vyb08}, comparable to the technique from \cite{Moe23}.

We begin by introducing weighted Wiener spaces and some of their properties in more detail in Section \ref{sec:wW_spaces}. After that we define some different \(s\)-numbers and related quantities such as the best \(m\)-term approximation in Section \ref{sec:s-numbers}. 
In Section \ref{sec:E_M} we develop some auxiliary results for the error of best approximation between Wiener spaces, before we finally have the first main results in Section \ref{sec:gel_ca}. Here we will give sharp bounds on the behavior of Gelfand numbers between Wiener spaces. In particular we show in Theorem \ref{thm:gel_th<eta<2} for \(0< \theta \leq \eta \leq 2\) with \(r > \Big(\frac{2 (\eta -\theta)(\theta -1)}{\eta (2 -\theta) \theta}\Big)_+\) (for \(\theta \leq 1\) this is zero), that it holds
    \begin{equation}
         c_m(\ca^{r}_\theta )_{\ca_\eta} \asymp m^{-(r  - \frac{1}{\eta} + \frac{1}{\theta})} \log(m)^{(d-1)r}.
    \end{equation}
We also show the same behavior in the reverse direction, see Theorem \ref{thm:gel_th>eta}. In Theorem \ref{thm:gel_2<th<eta} we cover the regime where \(2 \leq \theta \leq \eta \leq \infty\) and show that the rates are here independent of \(\theta\) and \(\eta\). More precisely for any \(r > 0\) we have
    \begin{equation}
         c_m(\ca^{r}_\theta)_{\ca_\eta} \asymp m^{-r} \log(m)^{(d-1)r}.
    \end{equation}

In Section \ref{sec:gel_Lq} we then use these results to construct bounds for the Gelfand numbers of \(\ca^r_\theta\) measured in \(L_q\). Here we distinguish between two different regimes. First we have \(1 \leq q \leq 2\) where the Gelfand numbers behave just as for \(q=2\). Indeed from Theorem \ref{thm:gel_L_q<2} we have for \(1 \leq \theta \leq 2\) and \(1 \leq q \leq 2\) with \(r > 1 - 1/\theta\) that it holds
    \begin{equation}
        c_m(\ca_\theta^r)_{L_q} \asymp m^{-(r -\frac{1}{2} + \frac{1}{\theta})} \log(m)^{(d-1) r}.
    \end{equation}
    The central result of this paper, Theorem \ref{thm:gel_L_q>2}, gives the rates for \(2 \leq q < \infty\) instead. Here we allow for \(0 <\theta \leq2\) but need a much stronger condition on the smoothness \(r > \frac{1}{2} - \frac{1}{q} + \big(1 - \frac{1}{\theta}\big)_+\). Under these assumptions it holds that
    \begin{equation}
      c_m(\ca_\theta^r)_{L_q} \asymp m^{-(r - 1 + \frac{1}{q} + \frac{1}{\theta})}\log(m)^{(d-1)(r - \frac{1}{2} +\frac{1}{q})}.
    \end{equation} 
The crux of this result is the gain in the logarithmic rates for \(q >2\). This is achieved by rewriting the \(L_q\)-norm akin to a Besov norm with dominating mixed smoothness and then applying the results for Gelfand numbers between Wiener spaces on each layer via the Hausdorff-Young inequality. For \(q = \infty\) the main rates can still be recovered but the logarithmic rates are simply \((d-1)r\).

In Section \ref{sec:samp_Lq} we can finally give lower bounds for the sampling numbers in \(L_q\) and compare them to the upper bounds from \cite{MSU26}. 
 For \(2 \leq q < \infty\), \(0 <\theta \leq 2\) and \(r > \frac{1}{2} - \frac{1}{q} +(1 - \frac{1}{\theta})_+\) we show
    \begin{equation}
         \varrho_m(\ca^r_\theta)_{L_q} \gtrsim m^{-(r - 1 + \frac{1}{q} + \frac{1}{\theta})}\log(m)^{(d-1)(r -  \frac{1}{2} + \frac{1}{q})}.
    \end{equation}
The main rates match the existing upper bounds, compare with \eqref{eq:intro_samp_upper}, and the gap in the non-dimension-dependent part of the logarithmic rates is to be expected. The gap \(\log(m)^{(d-1)(\frac{1}{2}- \frac{1}{q})}\) is significant however, and not something that could be easily bridged. Both the results for the best \(m\)-term approximation (in \(\ca\)) from \cite{MSU26} and the results for Gelfand numbers in \(L_q\) are sharp. Determining the precise (dimension-dependent) logarithmic rates for the sampling numbers of weighted Wiener spaces in \(L_q\) with \(2 < q < \infty\) remains an open problem.

{\bf Notation.} As usual $\N$ denotes the natural numbers, $\N_0=\N\cup\{0\}$,
$\zz$ denotes the integers,
$\re$ the real numbers, and $\C$ the complex numbers. The letter $d$ is always
reserved for the underlying dimension in $\re^d, \zz^d$ etc. For $a\in \re$ we denote $a_+ \coloneqq \max\{a,0\}$.
For $0<p\leq \infty$ and $\bx\in \R$ we denote $\Vert\bx\Vert_p = (\sum_{i=1}^d
|x_i|^p)^{1/p}$ with the
usual modification in the case $p=\infty$. By \(\bx \by \) we indicate the inner product. If $X$ and $Y$ are two (quasi-)normed spaces, the (quasi-)norm of an element $x$ in $X$ will be denoted by \(\Vert x\Vert_X\). We use the usual notations \(L_p\) for the Lebesgue (quasi-)Banach spaces and \(\ell_p\) for the sequence spaces, we abbreviate the \(p\) norm of \(x\) by \(\Vert x\Vert_p\).
The symbol $X \hookrightarrow Y$ indicates that the
identity operator is continuous and \(X'\) denotes the dual space of \(X\). For \(1 \leq p \leq \infty\) we denote by \(p'\) the conjugate exponent such that \(1 = \frac{1}{p} + \frac{1}{p'}\). For two sequences $a_n$ and $b_n$ we will write
$a_n \lesssim b_n$ if there exists a constant $c>0$, which may depend on \(d,r,s,\theta,\eta,q, \kappa\) but never on \(m,n,M\), and \(n_0 \in \N \) such that $a_n \leq c\,b_n$
for all $n \geq n_0$. We will write $a_n \asymp b_n$ if $a_n \lesssim b_n$ and $b_n
\lesssim a_n$.  
By \(\tor\) we indicate the set \( [0,1]\) where we identify the endpoints and denote by $f\in \mathcal{T}(Q)$ that \(f\) is a trigonometric polynomial with support on the frequencies in \(Q \subset \Z\). By \(B^r_{p,\theta}\) we denote the Besov spaces with dominating mixed smoothness.

\section{Weighted Wiener spaces}\label{sec:wW_spaces}

Weighted Wiener spaces are function spaces on the \(d\)-dimensional torus that are defined via a weighted \(\ell_\theta\) norm over their weighted Fourier coefficients. Their weights have a form similar to dominating mixed smoothness Besov spaces and have embeddings into these for certain ranges of parameters, see \cite[Section 5]{MSU26}. They can also be seen as the periodic counterpart to Barron classes, see \cite{Vo22}.

 \begin{defi}[Weighted Wiener space] \label{Wmix}
For \( r \geq 0\) and \(0 < \theta \leq \infty\) we define the (quasi-)norm 
 \begin{equation}
     \|f\|_{\ca^{r}_{\theta} } \coloneqq \Big( 
\sum_{\bk\in\Z} \prod_{i=1}^d (1+|k_i|)^{ r\theta}|\hat{f}(\bk)|^\theta \Big)^\frac{1}{\theta}. \end{equation}
The corresponding space
 \begin{equation}
\ca^{r}_{\theta}(\tor^d) \coloneqq \ca^{r}_{\theta} \coloneqq \Big\{ f \in L_1(\tor^d)~:~\Vert f\Vert_{\ca^{r}_{\theta} } < \infty \Big\}  
 \end{equation}
is called weighted (quasi-)Wiener space of functions with mixed smoothness.
In the case \(\theta = \infty\) the usual modification is made and the space is called Korobov space (see, e.g., \cite[Section 3.3]{DTU18}). The space \(\ca^0_1\) is simply the (original) Wiener space \(\ca\), and $\ca_{\theta} :=\ca^{0}_{\theta} $. The mixed Wiener spaces \(\ca_1^r \) with \(\theta = 1\) are often referred to as \(\ca^r_\mix\) in the literature.
 \end{defi}

Now let us introduce some of the usual tools that we will use during our proofs.

\begin{rem}\label{rem:embed}
While the space \(\ca^{r}_{\theta}\) is a function space, it can also be embedded into \(\ell_{\theta}\) quite naturally for all \(\theta > 0\) by the following  {isometric} operator 
 \begin{equation}\label{eq:Ar}
     A_r f= \bigg(\prod_{j=1}^d (1+|k_j|)^r \hat{f}(\bk)\bigg)_{\bk\in \zz^d}.
 \end{equation}
\end{rem}

We consider the product weights from the \(\ca\)-norm and introduce a notation for them:
\begin{equation}\label{w_k}
\omega_\bk \coloneqq \prod_{j = 1}^d (\vert k_j \vert +1) .
\end{equation}
The norm of our spaces \(\ca_\theta^r\) can now be rewritten as 
\begin{equation}
    \|f\|_{\ca^{r}_{\theta}} = \big( \sum_{\bk \in \Z} \omega_\bk^{ r \theta}\big\vert \hat{f}(\bk) \big\vert^\theta \big)^\frac{1}{\theta}.
\end{equation}

\begin{defi}[Uniform hyperbolic layers and blocks] \label{cross}
We set \(\hyp_0 = \{0\} \subset \Z \) for \(d \geq 1\) and define the uniform hyperbolic layer \(\hyp_j\) via the segmentation \(\zz = \bigcup_{n=-\infty}^\infty I_n\) 
with $I_0 = \{0\}$ and for $n\in \N$
 $$
 I_n = \big\{k \in \zz ~:~  2^{n-1} \leq  k < 2^{n}\big\}\,.
 $$
 This gives \(\#I_n = 2^{n-1}\).
 For \(n \in \zz\) with \(n< 0\) we extend this definition by setting \(I_n = - I_{-n}\).
For \( \bk \in \Z \) we write 
\begin{equation}\label{Ik}
    I_{\bk} = I_{k_1} \times \ldots \times I_{k_d},
\end{equation}
each of these blocks contains 
 \begin{equation}\label{Vol_Ik}
 \#I_{\bk} = 2^{\Vert\bk \Vert_1 -d }  
 \end{equation}
points. This also immediately induces a segmentation of \(\Z\) into cuboids, that is \(\Z = \bigcup_{\bk \in \Z} I_{\bk} \).

Now we can group blocks on the same hyperbolic layer together 
\begin{equation} \label{square}
 \hyp_j \coloneqq \bigsqcup_{\bk \in \Z, \, \Vert \bk \Vert_1 = j}  I_{\bk} ,
\end{equation}
where the contents of every \(I_{\bk}\) appear in exactly one \(\hyp_j, \, j \in \N_0\), meaning \(\Z = \bigsqcup_{j=0}^\infty \hyp_j\). 
Finally we may group all hyperbolic layers up to some level together and obtain a full hyperbolic cross
\begin{equation}
    \hyp^M \coloneqq \bigsqcup_{j = 0}^M \hyp_j.
\end{equation}

\end{defi}

Throughout this paper we will also need the number of points contained in a single hyperbolic layer \(\hyp_j\). We know that \(\hyp_j\) can be decomposed into a number of blocks \(I_{\bk}\) that each contain \(2^{j-d}\) points. We therefore have \(2^{j-d} 2^{d} = 2^j\) times the number of possibilities to distribute \(j \in \N\) to \(d\) dimensions. In total 
\begin{equation}\label{C_j}
C(j) \coloneqq \# \hyp_j = 2^j \binom{j+d-1}{j} \asymp 2^j j^{d-1}.
\end{equation}

We now define the \enquote{inverse} quantity to \(C(j)\), 
\begin{equation}\label{L}
L(n)= L = \big\lfloor n - (d-1) \log n \big\rfloor,
\end{equation}
this is chosen in such a way that
\begin{equation}\label{L_upper}
    C(L(n)) \asymp 2^n. 
\end{equation}

In particular there exists a constant \(c_d\) depending only on \(d\) such that \[C(L + c_d) \geq 2^{n+1}.\]

\bigskip

Another important property of Wiener spaces that we will use is that they can be directly related to \(L_p\) spaces by the Hausdorff-Young inequality.
\begin{prop}[Hausdorff-Young, {\cite[Satz I and II]{Haus23}}]\label{prop:Haus_young}
    Let  \(1 \leq q \leq 2\) and \(\frac{1}{q} + \frac{1}{q'} = 1\). Then it holds that 
    \begin{equation}
        \ca_{q} \hookrightarrow L_{q'}.
    \end{equation}
    In addition we have
        \begin{equation}
       L_q \hookrightarrow \ca_{q'}.
        \end{equation}
       In particular for any \(f \in L_1\)
       \begin{equation}
           \Vert f \Vert_{L_{q'}} \leq \Vert f \Vert_{\ca_q} 
       \end{equation}
       as well as
        \begin{equation}
       \Vert f\Vert_{ \ca_{q'}} \leq \Vert f\Vert_{L_q}.
        \end{equation}
\end{prop}

\section{Some \(s\)-numbers}\label{sec:s-numbers}

In this section we will introduce some \(s\)-numbers and widths, see, e.g., \cite{Pie87}. In particular we will consider the non-linear sampling widths and Gelfand numbers.

Let us first introduce the (non-linear) sampling widths. They give the worst case error of function recovery for some function class using standard information (point evaluations) and an arbitrary recovery algorithm.

\begin{defi}[Non-linear sampling widths] \label{defi_nl_samp}
Let $\mathcal{F}$ be a (quasi-)normed space of functions on $D$, where function evaluations are continuous, which is continuously embedded into a Banach space $Y$. Then the \(m\)-th (non-linear) sampling width is defined as
\begin{equation}  \label{defi_samp_eq}
\varrho_m(\mathcal{F})_Y \coloneqq \inf_{t_1, \ldots,  {t_{m}} \in D} \inf_{R\colon {\C^{m}} \to Y} \sup_{\|f\|_\mathcal{F} \leq 1} \big\Vert f - R \big(f(t_1), \ldots, f( {t_{m}})\big)\big\Vert_Y.
\end{equation}
\end{defi}

We now turn to the Gelfand numbers which allow for the same general recovery methods but also much more general information.

\begin{defi}[Gelfand numbers]\label{def:gel_num}
Let $X,Y$ be quasi-normed spaces and $T\colon X\to Y$ linear and continuous. The Gelfand numbers $c_m(T)$ for \(m \in \N_0\) are defined as 
\begin{equation}
c_m(T:X\to Y) \coloneqq \inf\Big\{\sup_{\Vert x\Vert_{X\cap M}\leq 1}\|Tx\|_Y\colon M\subset X \textrm{ closed linear,}\  \operatorname{codim} M\leq m \Big\}.
\end{equation}
The codimension of a subspace is defined as the dimension of the quotient space $X/M$.
\end{defi}

Since we want to use the Gelfand numbers to establish lower bounds on the sampling widths we also want to have a formulation of these as widths. We could formulate the corresponding width by plugging in the identity \(\id:X \to Y\) for the operator \(T\) but instead we will use a slightly different notation.

\begin{defi}[Gelfand widths]\label{def:gel_wid}
Let $\mathcal{F}$ be a (quasi-)normed space of functions continuously embedded into a Banach space $Y$. 
The $m$-th Gelfand width $c_m(\mathcal{F})_Y$ is defined as  
\begin{equation}
c_m(\mathcal{F})_Y \coloneqq \inf_{\lambda_1, \ldots, \lambda_{m} \in \mathcal{F}'} \inf_{R\colon\C^{m} \to Y} \sup_{\|f\|_\mathcal{F} \leq 1} \big\Vert f - R \big(\lambda_1(f), \ldots, \lambda_{m}(f)\big)\big\Vert_Y.
\end{equation}
\end{defi}

We now want to justify these definitions by showing that they differ by just a constant factor. This was done in e.g., \cite[Theorem 10.4]{FoRa13} with a focus on adaptivity.
\begin{lem}
 If $\mathcal{F}$ is a normed space (of functions) and $\id:\mathcal{F} \to Y$, where $Y$ is also a normed space, then 
$$
 \frac{1}{2}c_m(\mathcal{F})_Y \leq c_m(\id) \leq c_m(\mathcal{F})_Y, 
$$ 
where the latter quantities are the Gelfand widths from Definition \ref{def:gel_wid} above. In case of quasi-Banach spaces this identity does not hold as pointed out after \cite[Def.\ 2.1]{DiUl18}. However, we still have 
$$
c_m(\id) \leq c_m(\mathcal{F})_Y, 
$$ 
even when \(X\) is only quasi-Banach.
\end{lem}

As already implied all of the quantities so far do satisfy the definition of an \(s\)-width or \(s\)-number, see, e.g., \cite{Pie87}.

\begin{defi}[s-function]\label{def:sfunc}
Let \(X, Z\) be quasi-Banach spaces and $Y$ be a $p$-Banach space for \(0 <p \leq 1\), let $S,T\in\mathcal{L}(X,Y)$ and $R\in \mathcal{L}(Y,Z)$. A mapping $s\colon T \mapsto (s_m(T))_{m=0}^{\infty}$ with the following properties
\begin{description}
  \item[(S1)] $\|T\|_{\mathcal{L}(X,Y)}=s_0(T)\geq s_1(T)\geq\ldots\geq 0$
  \item[(S2)] for all $ m_1,m_2\in\N_0 $ it holds
  $$s_{m_1+m_2}(R\circ S)\leq s_{m_1}(R)s_{m_2}(S)$$
  \item[(S3)] for all $ m_1,m_2\in\N_0 $ it holds $$s_{m_1+m_2}^p(S+T)\leq s_{m_1}^p(S)+s_{m_2}^p(T)$$
  \item[(S4)] If $\rank T \leq m$ then $s_m(T)=0$ and $s_{m}(\id\colon\ell_2^{m+1}\to \ell_2^{m+1})=1$
\end{description}
is called s-function and \(s_m(T)\) is called the \(m\)-th s-number of the operator \(T\).
\end{defi}

\begin{prop} \label{prop:all_s_num}
    Gelfand numbers are \(s\)-numbers.
\end{prop}
\begin{proof}
    See, e.g., \cite[Theorem 11.5.5]{Pie78} and for the quasi-Banach case see also \cite[Section 2]{DiUl18}.
\end{proof}

As mentioned already in the beginning of this section the Gelfand widths provide a lower bound on the sampling widths. This is almost immediately clear for the formulations of the widths used here but was also shown in e.g., \cite[Lemma B.1]{JUV23}.

\begin{lem}\label{lem:gel_samp}
    Let \(X\) be a quasi-Banach space of functions which is continuously embedded into \(C(\tor^d)\) and \(Y\) be another Banach space such that \(C(\tor^d) \hookrightarrow Y\). Then it holds for all \(m \in \N\) that,
    \begin{equation}
        c_m(X)_Y \leq \varrho_m(X)_Y.
    \end{equation}
    In particular we have that this result is applicable for \(\ca^r_\theta \hookrightarrow \ca \hookrightarrow C(\tor^d) \hookrightarrow L_q\), for \(r > (1 - 1/\theta)_+\).
\end{lem}

\bigskip

Finally we want to define some trigonometric approximation characteristics. These do not satisfy the definition of an \(s\)-width.

The best \(m\)-term approximation will not be studied in this paper. Results for this quantity in weighted Wiener spaces were already shown in \cite{MSU26} based on results for unweighted Wiener spaces, see, e.g., \cite{MSU25}.

\begin{defi}[Best trigonometric \(m\)-term  approximation width]\label{def:bm}
Let \(X,Y\) be quasi-Banach function spaces on \(\tor^d\)  and \(f \in X\).
The best \(m\)-term  (trigonometric) approximation (width) of \(f\) is defined as 
\begin{equation}
    \sigma_m(f)_Y \coloneqq \inf_{s \in \Sigma_{m}}  \|f - s\|_Y, 
\end{equation} 
where \(\Sigma_m\) is the set of all trigonometric polynomials of the form 
 \[\sum_{\bk \in \zz^d} a_\bk \exp(2 \pi \mathrm{i}\bk\bx) \] 
with at most \(m\) non-zero coefficients \(a_\bk\).

The best $m$-term approximation width of a function class \(\cF \hookrightarrow Y\) is defined accordingly, as
\[ \sigma_m(\cF)_Y \coloneqq \sup_{ \Vert f \Vert_\cF \leq 1} \sigma_m(f)_Y.\]
\end{defi}

\begin{defi}[Error of best trigonometric approximation]\label{defi:E_M}
Let \(Q \subset \zz^d\). The error of best trigonometric approximation of a function space \( \cF \hookrightarrow Y \hookrightarrow L_1\) is defined as
\begin{equation} 
E_{Q}(\cF)_{Y} \hspace{-2pt}
 \coloneqq \hspace{-2pt}  {\sup_{\Vert f \Vert_{\cF} \leq 1} } 
  {E_{Q}\big(f\big)_{Y}
 \coloneqq }\hspace{-2pt}  {\sup_{\Vert f \Vert_{\cF} \leq 1} }\inf_{g \in \mathcal{T}(Q)} \Vert f-g \Vert_{Y}\,.
\end{equation}
\end{defi}

\section{Error of best Trigonometric Approximation}\label{sec:E_M}

This section is dedicated to observing the behavior of the error of best trigonometric approximation, see Definition \ref{defi:E_M}. This will be useful later on, when we study Gelfand numbers as it is a natural upper bound for them. The error of best trigonometric approximation, measured in \(L_\infty\), also appears in upper bounds on the sampling numbers from \cite{JUV23, Kri23, MSU26}, see Proposition \ref{prop:samp_A1}.
The following result generalizes \cite[Lemma 4.6]{MSU26} in multiple directions but uses a similar basic strategy. It covers a greater range of parameters and gives lower bounds in the setting where the target space is \(\ca_\eta\).
\begin{thm}\label{thm:E_M}
    For \(\eta, \theta \in (0,\infty]\), \(s\geq 0\) and \( r > \big(\frac{1}{\eta} - \frac{1}{\theta}\big)_+\) it holds
    \begin{equation}
        E_{[-M,M]^d}(\ca^{r+s}_\theta)_{\ca^s_\eta} = E_{[-M,M]^d}(\ca^{r}_\theta)_{\ca_\eta} \asymp 
        \begin{cases}
            M^{-r} & \theta \leq \eta \\
            M^{-(r -\frac{1}{\eta}+ \frac{1}{\theta})} & \theta > \eta.
        \end{cases}
    \end{equation}
\end{thm}
\begin{proof}
We start by considering the setting where \(\theta \leq \eta\). It is immediately obvious that the function \(f(\bx) = (M+2)^{-r} \exp(2\pi \mathrm{i} (M+1) x_1)\) realizes the lower bound. The upper bound is trivial by the embedding \(\ell_\theta \hookrightarrow \ell_\eta\). Therefore we have \(E_{[-M,M]^d}(\ca^r_\theta)_{\ca_\eta} \asymp M^{-r}\). Now let us turn to the more involved regime where \(\theta > \eta\), here we know that  \(r \frac{\theta\eta}{\theta-\eta} > 1\). Now apply Hölder's inequality with \(\frac{\theta}{\eta} > 1\),
    \begin{align}
        \begin{split}
             E_{[-M,M]^d}(\ca^r_\theta)_{\ca_\eta} & = \sup_{\Vert f\Vert_{\ca^r_\theta} \leq 1} \Big(\sum_{\bk \in \Z \setminus [-M,M]^d} \omega^{-r\eta}_\bk \omega^{r\eta}_\bk \vert \hat{f}(\bk)\vert^\eta \Big)^{1/\eta} \\
             & \leq \sup_{\Vert f\Vert_{\ca^r_\theta} \leq 1}\bigg( \Big( \sum_{\bk \in \Z \setminus [-M,M]^d} \hspace{-3pt} \omega^{-r\frac{\theta\eta}{\theta-\eta}}_\bk \Big)^{\frac{\theta-\eta}{\theta}} \Big( \sum_{\bk \in \Z \setminus [-M,M]^d} \hspace{-3pt} (\omega^{r\eta}_\bk \vert \hat{f}(\bk) \vert^{\eta} )^{\frac{\theta}{\eta}} \Big)^{\frac{\eta}{\theta}} \bigg)^{1/\eta} \\
             & = \sup_{\Vert f\Vert_{\ca^r_\theta} \leq 1} \Big( \sum_{\bk \in \Z \setminus [-M,M]^d} \hspace{-3pt} \omega^{-r\frac{\theta\eta}{\theta-\eta}}_\bk \Big)^{\frac{\theta-\eta}{\theta\eta}} \Big( \sum_{\bk \in \Z \setminus [-M,M]^d} \hspace{-3pt} \omega^{r\theta}_\bk \vert \hat{f}(\bk) \vert^{\theta} \Big)^{\frac{1}{\theta}}  \\
             & \leq  \Big( \sum_{\bk \in \Z \setminus [-M,M]^d} \omega^{-r{\frac{\theta\eta}{\theta-\eta}}}_\bk \Big)^{\frac{\theta-\eta}{\theta\eta}} \\[-5pt]
            & \leq\Big( d \sum_{\bk \in \zz^{d-1}} \sum_{l \in \zz \setminus [-M,M]} (1+ \vert l\vert)^{-r\frac{\theta\eta}{\theta-\eta}} \prod_{j=1}^{d-1} (1+\vert k_j \vert)^{-r\frac{\theta\eta}{\theta-\eta}} \Big)^{\frac{\theta-\eta}{\theta\eta}} \\
            & \leq \Big( 2d \sum_{l = M +1}^\infty (1+ \vert l\vert)^{-r\frac{\theta\eta}{\theta-\eta}} \sum_{\bk \in \zz^{d-1}}  \prod_{j=1}^{d-1} (1+\vert k_j \vert)^{-r\frac{\theta\eta}{\theta-\eta}} \Big)^{\frac{\theta-\eta}{\theta\eta}} \\
            &\lesssim \Big(\sum_{l = M +1}^\infty (1+ \vert l\vert)^{-r\frac{\theta\eta}{\theta-\eta}} \Big)^{\frac{\theta-\eta}{\theta\eta}} \\
            & \lesssim M^{\frac{1}{\eta}- \frac{1}{\theta}-r}. \\[-1pt]
        \end{split}
    \end{align}   Where in the second to last line we used the fact that \( \sum_{\bk \in \zz^{d}}  \prod_{j=1}^{d} (1+\vert k_j \vert)^{-r{\frac{\theta\eta}{\theta-\eta}}} \leq C^d\) for some \(C\) (this can easily be checked by induction over \(d\)). The inequality in the last line holds by comparing to the integral.

To give a lower bound we provide a fooling function from \(\ca_\theta^r(\tor)\). Using a function on the one-dimensional torus is sufficient since the behavior of the error is (almost) independent of the underlying dimension. We set
\begin{equation}
        f(x) = \sum_{k = M+1 }^{2M} k^{-r} M^{-1/\theta} e^{2 \pi\mathrm{i} k x}.
    \end{equation}
    Since we have \(\supp(\hat f) \cap [-M,M] = \emptyset\) we know that the best possible approximation from \(\mathcal{T}([-M,M])\) in \(\ca_\eta\) is simply doing nothing, i.e., \(g = 0\) in Definition \ref{defi:E_M}. Therefore we have 
     \begin{align}
        \begin{split}
         E_{[-M,M]^d}(\ca^r_\theta)_{\ca_\eta} & \geq  E_{[-M,M]}(f)_{\ca_\eta} = \Vert f \Vert_{\ca_\eta} \\
         & = \big(\sum_{k = M +1}^{2M} k^{-r\eta} M^{-\eta/\theta} \big)^{1/\eta}\\
         & \geq (2M)^{-r} M^{\frac{1}{\eta}-\frac{1}{\theta}}\\
         & \gtrsim M^{\frac{1}{\eta}-\frac{1}{\theta}-r}.\\
        \end{split}
    \end{align} 
    The bounds still hold for \(\theta = \infty\) or \( \eta = \infty\) with the usual modifications.
\end{proof}
We are also interested in the situation where \(Q = \hyp^M\) is a hyperbolic cross.
\begin{thm}\label{thm:E_M_hyp}
        For \(\eta, \theta \in (0,\infty]\), \(s\geq 0\) and \( r > \big(\frac{1}{\eta} - \frac{1}{\theta}\big)_+\) it holds
    \begin{equation}
        E_{\hyp^M}(\ca^{r+s}_\theta)_{\ca^s_\eta} = E_{\hyp^M}(\ca^{r}_\theta)_{\ca_\eta} \asymp 
        \begin{cases}
            2^{-Mr}, & \theta \leq \eta \\
            2^{-M(r -\frac{1}{\eta}+ \frac{1}{\theta})} M^{(\frac{1}{\eta}- \frac{1}{\theta})(d-1)}, & \theta > \eta.
        \end{cases}
    \end{equation}
\end{thm}
\begin{proof}
The first case is again trivial. We now consider \(\theta > \eta\), because our target space is of \(\ca\)-type we again have that the error of best trigonometric approximation coincides with the projection error. We use that by Hölder's inequality the \(\ell_\eta\)-norm of a sequence with at most \(C(j)\) non-zero entries can be estimated from above by its \(\ell_\theta\)-norm times \(C(j)^{1/\eta-1/\theta}\) (see, e.g., \cite[Lemma 4.3]{MSU26}),
 
\begin{align} \label{eq:E_M_hyp_block}
    \begin{split}
        E_{\hyp^M}(\ca^{r}_\theta)_{\ca_\eta}^\eta & = \sup_{\Vert f \Vert_{\ca_\theta^r} \leq 1} \sum_{j = M+1}^\infty \Big(\sum_{\bk \in \hyp_j} \vert \hat{f}(\bk) \vert^\eta \Big)^{\eta/\eta} \\
        & \leq \sup_{\Vert f \Vert_{\ca_\theta^r} \leq 1} \sum_{j = M+1}^\infty C(j)^{1-\eta/\theta}\Big(\sum_{\bk \in \hyp_j} \vert \hat{f}(\bk) \vert^\theta \Big)^{\eta/\theta} \\
        & \lesssim \sup_{\Vert f \Vert_{\ca_\theta^r} \leq 1} \sum_{j = M+1}^\infty C(j)^{1-\eta/\theta} 2^{-\eta r j} \Big(  \sum_{\bk \in \hyp_j} \omega_{\bk}^{\theta r} \vert \hat{f}(\bk) \vert^\theta \Big)^{\eta/\theta} \\
        & \asymp  \sup_{\Vert f \Vert_{\ca_\theta^r} \leq 1} \sum_{j = M+1}^\infty (2^j j^{(d-1)})^{1-\eta/\theta} 2^{-\eta rj} \Big(  \sum_{\bk \in \hyp_j} \omega_{\bk}^{\theta r} \vert \hat{f}(\bk) \vert^\theta \Big)^{\eta/\theta} \\
        & \leq \sum_{j = M+1}^\infty 2^{j\eta(1/\eta-1/\theta-r)} j^{\eta(1/\eta-1/\theta)(d-1)}.
    \end{split}
\end{align}
  Since the sum is falling geometrically we can estimate it from above by some constant times the first summand (see, e.g., \cite[Lemma 4.2]{MSU26}). Then we get after taking the \(\eta\)-root
    \begin{equation}
E_{\hyp^M}\big(\ca_\theta^{r}\big)_{\ca_\eta}  \lesssim 2^{M(1/\eta-1/\theta-r)} M^{(1/\eta-1/\theta)(d-1)}.
    \end{equation}

     To construct a lower bound we need a fooling function.
 \begin{equation}
        f(\bx) = \sum_{\bk \in \hyp_{M+1}} \omega_\bk^{-r} C(M+1)^{-1/\theta} e^{2 \pi \mathrm{i} \bk \bx}.
 \end{equation}
    
    From \eqref{C_j} it is clear that \(\Vert f\Vert_{\ca_\theta^r}\leq 1\), so we only need to consider the error of best trigonometric approximation in our target space \(\ca_\eta\). We continue by noting 
    
    \begin{align} 
    \begin{split}
E_{\hyp^M}(\ca^{r}_\theta)_{\ca_\eta} &\geq E_{\hyp^M}(f)_{\ca_\eta}\\
& =  \Big( \sum_{\bk \in \hyp_{M+1}} \omega_\bk^{-r\eta} C(M+1)^{-\eta/\theta}\Big)^{1/\eta} \\
& \gtrsim \Big( C(M+1) \,2^{- Mr\eta} C(M+1)^{-\eta/\theta} \Big)^{1/\eta} \\
&\gtrsim 2^{-Mr} C(M)^{1/\eta-1/\theta} \\
&\gtrsim  2^{-M(r - \frac{1}{\eta}+ \frac{1}{\theta})} M^{(\frac{1}{\eta}- \frac{1}{\theta})(d-1)}.
    \end{split}
    \end{align}
    Where in the last line we used the definition of \(C(M)\), see \eqref{C_j}.
\end{proof}

\section{Gelfand numbers in Wiener spaces}\label{sec:gel_ca}

In this section we want to study the Gelfand numbers of embeddings \(\ca^r_\theta \to \ca^s_\eta\).
Bounds for \(\ca^r_\theta \hookrightarrow L_2 = \ca_2\) were studied in \cite[Theorem 3.2 and Corollary 3.3]{Moe23} based on bounds for Gelfand numbers from \cite[Theorem 1.1]{FPRU10}. In that work only the results for \(c_n(\ell^N_1 )_{\ell^N_2}\) from \cite{FPRU10} were used. However, since we are now interested in embeddings into all \(\ca_\eta\) spaces we can easily extend these results.
We start just as in \cite{Moe23} by considering that we can study Gelfand numbers of a diagonal operator between \(\ell_p\)-spaces instead of the identity between \(\ca_p\) spaces.

\begin{thm} \label{thm:gel_emb}
    Let \(\ca_\theta^r \hookrightarrow \ca_\eta^s\), that is \( \theta \leq \eta\) and \(r \geq s \geq 0\) or \( \theta \geq \eta\) and \(r > s + \frac{1}{\eta}- \frac{1}{\theta}\).
Then the following diagrams commute and we have for \(m \in \n\) that
\begin{equation}
    c_m(\id:\ca^{r}_\theta \to \ca^s_\eta) = c_m(D_{r-s}:\ell_\theta \to \ell_\eta).
\end{equation}

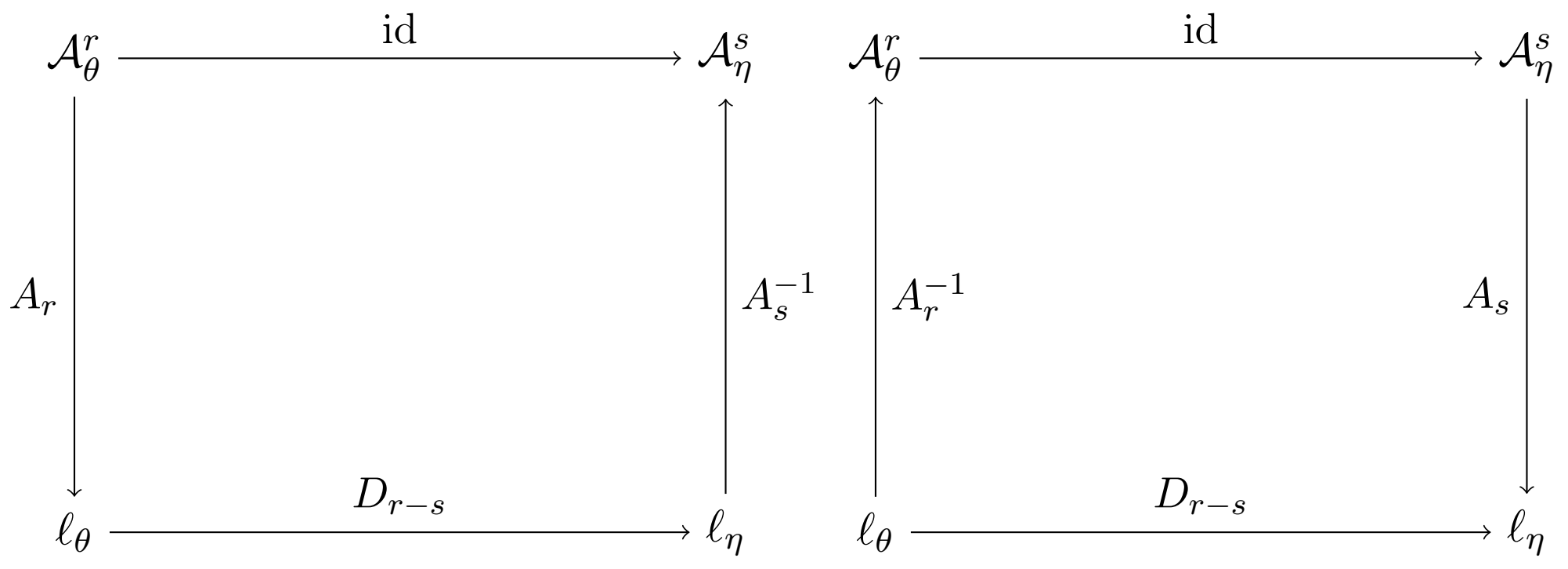
\begin{figure}[ht]
\begin{center} 
  \begin{tikzpicture} 
  \def\a{2.75} \def\b{2}
  \path
  (-\a,\b) node (A) {$\ca_\theta^r$} 
  (-\a,-\b) node (B) {$\ell_\theta$}      
  (\a,\b) node (C) {$\ca_\eta^s$}
  (\a,-\b) node (D) {$\ell_\eta$};
  \begin{scope}[nodes={midway,scale=.95}]
    \draw[->] (A)--(B) node [left] {$A_r$};
    \draw[->] (D)--(C) node [right] {\(A^{-1}_s\)};
    \draw[->] (A)--(C) node [above] {$\id$};
    \draw[->] (B)--(D) node [above] {\(D_{r-s}\)};
  \end{scope}
\end{tikzpicture}%
  \begin{tikzpicture} 
  \def\a{2.75} \def\b{2}
  \path
  (-\a,\b) node (A) {$\ca_\theta^r$} 
  (-\a,-\b) node (B) {$\ell_\theta$}      
  (\a,\b) node (C) {$\ca_\eta^s$}
  (\a,-\b) node (D) {$\ell_\eta$};
  \begin{scope}[nodes={midway,scale=.95}]
    \draw[->] (B)--(A) node [right] {$A_r^{-1}$};
    \draw[->] (C)--(D) node [left] {\(A_s\)};
    \draw[->] (A)--(C) node [above] {$\id$};
    \draw[->] (B)--(D) node [above] {\(D_{r-s}\)};
  \end{scope}
\end{tikzpicture}%
\end{center}
    \caption{Commutative diagrams between Wiener and sequence spaces}
     \label{fig:com_Gel}
\end{figure}
Where for \(t \in \mathbb{R}\) we define \(D_t\) as 
\begin{equation}
    D_t(x_\bk)_{\bk\in\Z}= \big(\omega_\bk^{-t} x_\bk\big)_{\bk\in\Z},
\end{equation}
 and \(A_r\), \(A_s^{-1}\) are the operators from Remark \ref{rem:embed}.
\end{thm}

\begin{proof}
    We start by checking that the first diagram in Figure \ref{fig:com_Gel} commutes. Let \(f \in \ca^r_\theta\), then
    \begin{align}
        \begin{split}
            A^{-1}_s \circ D_{r-s} \circ A_r f & =  A^{-1}_s \circ D_{r-s} \big(\omega_\bk^r \hat{f}(\bk)\big)_{\bk\in \zz^d} \\
            & = A^{-1}_s \big(\omega_\bk^s \hat{f}(\bk)\big)_{\bk\in \zz^d} \\
            & =  \sum_{\bk\in \zz^d} \hat{f}(\bk) \exp(2 \pi \mathrm{i} \bk \cdot) \\
            & = f.
        \end{split}
    \end{align}
    It is now trivial to check that the second diagram in Figure \ref{fig:com_Gel} also commutes.
We now use the fact that Gelfand numbers are s-numbers, see Proposition \ref{prop:all_s_num} and apply (S2) and (S1) from Definition \ref{def:sfunc},
\begin{align}
    \begin{split}
        c_m(\id:\ca^r_\theta \to \ca^s_\eta) & \leq c_0(A_r: \ca^r_\theta \to \ell_\theta) c_m(D_{r-s}:\ell_\theta \to \ell_\eta) c_0(A^{-1}_s: \ell_\eta \to \ca^s_\eta) \\
        & = \Vert A_r \Vert_{\ca^r_\theta \to \ell_\theta} c_m(D_{r-s}:\ell_\theta \to \ell_\eta) \Vert A^{-1}_s \Vert_{\ell_\eta \to \ca^s_\eta} \\
        & = c_m(D_{r-s}:\ell_\theta \to \ell_\eta).
    \end{split}
\end{align}
We can now do the same calculation for the second diagram and get
\begin{align}
    \begin{split}
        c_m(D_{r-s}:\ell_\theta \to \ell_\eta) & \leq c_0(A_r^{-1}: \ell_\theta \to \ca^r_\theta) c_m(\id:\ca^r_\theta \to \ca^s_\eta) c_0(A_s:\ca^s_\eta \to \ell_\eta) \\
        & = \Vert A_r^{-1}\Vert_{\ell_\theta \to \ca^r_\theta} c_m(\id:\ca^r_\theta \to \ca^s_\eta) \Vert A_s \Vert_{\ca^s_\eta \to \ell_\eta} \\
        & = c_m(\id:\ca^r_\theta \to \ca^s_\eta).
    \end{split}
\end{align}
\end{proof}

We will now state a few very useful lemmata from the literature that we will use to show bounds for Gelfand numbers between Wiener spaces. All of these are in the setting of sequence spaces and stated for the identity. We will use our typical block decomposition method to apply them to the weighted regime. All of these are stated for cut-off \(\ell_p\) spaces.

\begin{defi}
For \(N \in \N\) and \(0 < p \leq \infty\) we define the cut-off sequence spaces 
\[\ell_p^N \coloneqq\Big\{(x_n)_{n=1}^{N}\Big|\|\bx\|_{\ell_p^N}\coloneqq\Big(\sum_{n=1}^{N}|x_n|^p\Big)^{\frac{1}{p}}<\infty\Big\},\]
with the usual modifications for \(p = \infty\).
\end{defi}

We begin with a result where the target norm is larger than the origin norm. 
\begin{lem}[{\cite[Lemma 4.8]{Vyb08}}]\label{lem:gel_rev}
Let \(0 < q \leq p \leq \infty\), then it holds for all \(N,m \in \N\) with \(N > m\) that
\begin{equation}
    c_m(\ell_p^N)_{\ell_q^N} = (N - m)^{1/q-1/p}.
\end{equation}
\end{lem}

Next up we have the classical regime including from \(\ell_1\) to \(\ell_2\), the upper bound was originally shown in \cite[Lemma 4.11]{Vyb08} and is based on results from \cite{Kas77}.
\begin{lem}[{\cite[Theorem 1.1]{FPRU10}}]\label{lem:gel_rates_ell} 
    Let \(0 < p \leq 1\) and \(p < q \leq 2\) then it holds for all \(N,m \in \N\) with \(N > m\) that
    \begin{equation}
        c_m(\ell_p^N)_{\ell_q^N} \asymp \Big( \frac{\log(e N/m)}{m} \Big)^{1/p - 1/q}.
    \end{equation}
\end{lem}
The regime where \(1 < p < q\) is instead governed by the following lemma from \cite[Lemma 4.7]{Vyb08} and \cite[Korollar 7.6]{Buch99}, based on results for Kolmogorov numbers \cite[\begin{otherlanguage*}{russian}Теорема 1\end{otherlanguage*}]{Glu83}. 

\begin{lem}\label{lem:gel_p>1}
    Let \(1 < p < q \leq 2\) then it holds for all \(N,m \in \N\) with \(N > m\) that
    \begin{equation}
        c_m(\ell_p^N)_{\ell_q^N} \asymp \min\big\{1, N^{1-1/p} m^{-1/2}\big\}^{\frac{1/p-1/q}{1/p -1/2}}.
    \end{equation}
\end{lem}

 Finally the regime where \(2 \leq q\) behaves very differently and was studied in \cite[Lemma 4.7]{Vyb08} and \cite[Korollar 7.6]{Buch99}.
\begin{lem}\label{lem:gel_1<p<2}
    Let \(1 < p \leq 2 \leq q \leq \infty\) then it holds for all \(N,m \in \N\) with \(N > m\) that
    \begin{equation}
        c_m(\ell_p^N)_{\ell_q^N} \asymp \max\big\{N^{1/q-1/p},\min\big\{1, N^{1-1/p} m^{-1/2}\big\} \sqrt{1- m/N}\big\} .
    \end{equation}
    In addition for \(p = 1\) it holds
    \begin{equation}
         c_m(\ell_1^N)_{\ell_q^N} \gtrsim \max\big\{N^{1/q-1}, m^{-1/2}\sqrt{ 1-m/N}\big\}.
    \end{equation}
\end{lem}

\begin{thm} \label{thm:gel_th>eta}
    Let \(s \geq 0\), \(0 < \eta \leq \theta \leq \infty\) and \(r > \frac{1}{\eta}- \frac{1}{\theta}\) then it holds 
    \begin{equation}
         c_m(\id:\ca^{r+s}_\theta \to \ca^s_\eta) \asymp m^{-(r  - \frac{1}{\eta} + \frac{1}{\theta})} \log(m)^{(d-1)r}.
    \end{equation}
\end{thm}
\begin{proof}
    We start by employing Theorem \ref{thm:gel_emb} and can now study \(c_m(D_{r}:\ell_\theta \to \ell_\eta)\) to get the desired rates. We follow the proof idea from \cite[Theorem 3.2]{Moe23} but modify it to suit the more general setting. 
    We start with our standard decomposition onto the layers \(\hyp_j\). We define the restriction of \(D_{r}\) to \(\hyp_j\) as 
    \begin{equation}
            \Delta_j(x_\bk)  \coloneqq \begin{cases}
                0 & \bk \notin \hyp_j \\
                \omega^{-r}_\bk x_\bk & \bk \in \hyp_j,
            \end{cases}
    \end{equation}
    for all \(\bk \in \Z\).
    Now we apply (S3) and write for \(p = \min\{1, \eta\}\),
    \begin{equation} \label{eq:decomp_gel_2}
        c^p_m(D_{r}) \lesssim \sum_{j=0}^{L} c^p_{m_j}(\Delta_j) + c^p_{0}\big(\sum_{j=L+1}^\infty\Delta_j\big).
    \end{equation}
   Let always \(m = 2^n\) and \(L\) be as in \eqref{L}. We set \(m_j = C(j)\) for \(j = 0, \ldots, L\). In this domain we therefore get \(c_{m_j}(\Delta_j) = 0\) by (S4).
   Then we have from \eqref{L_upper} that \(\sum_{j=0}^{L} C(j)\leq 2 C(L) \lesssim m\).
All that is left is to investigate the setting where \(j \geq L+1\). Here we set \(m_j = 0\) and by (S1) observe how the norm behaves,
\begin{align} 
    \begin{split}
        c_{0}\big(\sum_{j=L+1}^\infty\Delta_j\big) & = \Big\Vert  \sum_{j=L+1}^\infty\Delta_j \Big\Vert_{\ell_\theta \to \ell_\eta}\\
        & = \sup_{\Vert (x_\bk)_{\bk \in \Z} \Vert_{\ell_\theta} \leq 1} \Big(\sum_{\bk \in \Z \setminus \hyp^L} \vert\omega_\bk^{-r} x_\bk\vert^\eta \Big)^{1/\eta}\\
        & = \sup_{\Vert f \Vert_{\ca_\theta^{r}} \leq 1} \Vert f - f_{\hyp^L}\Vert_{\ca_\eta} \\
        & = E_{\hyp^L}(\ca^{r}_\theta)_{\ca_\eta}. \\
    \end{split}
\end{align}
    Where in the last line we used that the error of best trigonometric approximation and the trigonometric projection error agree here, since we measure in a \(\ca\)-type space.
    Now we can apply Theorem \ref{thm:E_M_hyp} with \(M = L\) and obtain
    \begin{align} \label{eq:delta_EL}
        \begin{split}
            c_{0}\big(\sum_{j=L+1}^\infty\Delta_j\big) & \asymp  2^{-L(r -\frac{1}{\eta}+\frac{1}{\theta})} L^{(\frac{1}{\eta}- \frac{1}{\theta})(d-1)} \\
            & \asymp 2^{-(n - (d-1) \log n)(r -\frac{1}{\eta}+\frac{1}{\theta})} (n - (d-1) \log n)^{(\frac{1}{\eta}- \frac{1}{\theta})(d-1)} \\
            & \leq 2^{-n(r -\frac{1}{\eta}+\frac{1}{\theta})} n^{(d-1)r}\\
            & =  m^{-(r -\frac{1}{\eta}+\frac{1}{\theta})} \log(m)^{(d-1)r}.
        \end{split}
    \end{align}
    To show lower bounds we will restrict to  \(\Delta_{L+c_d}\). For this layer we have that the number of points on it is at least \(C(L+c_d) \geq 2 m\). Now we can use (S2) from Definition \ref{def:sfunc} and apply Lemma \ref{lem:gel_rev} to obtain
    \begin{align}
        \begin{split}
            c_m(\id:\ca^{r+s}_\theta \to \ca^s_\eta) &\geq c_m(\Delta_{L+c_d}: \ell_\theta \to \ell_\eta) \\
            & \gtrsim 2^{-r(L+c_d)} c_m(\id: \ell_\theta^{C(L+c_d)} \to \ell_\eta^{C(L+c_d)}) \\
            & \gtrsim 2^{-rL} \big(C(L+c_d) -m\big)^{1/\eta-1/\theta} \\
            & \gtrsim 2^{-nr + \log(n)(d-1)r} (2m-m)^{1/\eta-1/\theta} \\
            &  = m^{-(r- 1/\eta + 1/\theta)} \log(m)^{(d-1)r}.
        \end{split}
    \end{align}
    This finishes the proof after arguing with monotonicity.
\end{proof}
The regime where \(\theta < \eta\) is a little bit more involved and follows the proof from \cite[Theorem 3.2 and Corollary 3.3]{Moe23} more closely. The main difference to the previous result is that \(E_{\hyp^L}\) does not admit the correct rate for the Gelfand numbers anymore and therefore we need to decompose the problem into three parts now, instead of into two like in \eqref{eq:decomp_gel_2}. To treat the newly appearing middle sum we will employ the above results for Gelfand numbers between cut-off \(\ell_p\) spaces.

\begin{thm} \label{thm:gel_th<eta<2}
    Let \(s \geq 0\) and \(0< \theta \leq \eta \leq 2\) with \(r > \Big(\frac{2 (\eta -\theta)(\theta -1)}{\eta (2 -\theta) \theta}\Big)_+\), then it holds 
    \begin{equation}
         c_m(\id:\ca^{r+s}_\theta \to \ca^s_\eta) \asymp m^{-(r  - \frac{1}{\eta} + \frac{1}{\theta})} \log(m)^{(d-1)r}.
    \end{equation}
\end{thm}
\begin{proof}
    We start just as in Theorem \ref{thm:gel_th>eta} by Theorem \ref{thm:gel_emb} however we now use the decomposition 
    \begin{equation} \label{eq:decomp_gel_3}
        c^p_m(D_{r}) \lesssim \sum_{j=0}^L c^p_{m_j}(\Delta_j) +  \sum_{j=L+1}^K c^p_{m_j}(\Delta_j)+  c^p_{0}\big(\sum_{j=K+1}^\infty\Delta_j\big),
    \end{equation}
    where again \(p = \min\{1, \eta\}\). By (S3) this inequality holds true under the condition that \(\sum_{j = 0}^K m_j \lesssim m\). Let as always \(m = 2^n\) and \(L\) be as in \eqref{L}. We again use the trivial approximation \(m_j = C(j)\) for \(j = 0, \ldots, L\). 
 Now we turn to \(j = L+1,\ldots K\). For \(\theta \leq 1\) we set \(m_j = \lceil 2^j 2^{(L-j)\kappa}j^{d-1} \rceil\) and apply Lemma \ref{lem:gel_rates_ell},
\begin{align} 
\begin{split}
\sum_{j=L+1}^K c_{m_j}^p(\Delta_j\colon \ell_\theta\to\ell_\eta) &\lesssim \sum_{j=L+1}^K c^p_{m_j}(\id\colon \ell_\theta^{C(j)}\to\ell_\eta^{C(j)})2^{-j r p} \\
&\lesssim\sum_{j=L+1}^K\Bigg(\frac{\log\Big(\frac{2^{j}j^{d-1}}{2^{j}2^{(L-j)\kappa}j^{d-1}}\Big)}{2^{j}2^{(L-j)\kappa}j^{d-1}}\Bigg)^{(\frac{1}{\theta}-\frac{1}{\eta}) p} 2^{-j r p}\\
& = \sum_{j=L+1}^K[(j-L)\kappa 2^{(j-L)\kappa}2^{-j}j^{1-d}]^{(\frac{1}{\theta}-\frac{1}{\eta}) p}2^{-j r p}\\
&\asymp 2^{-L\kappa(\frac{1}{\theta}-\frac{1}{\eta})p} \hspace{-2pt}\sum_{j=L+1}^K \hspace{-2pt} \Big[(j-L) 2^{-j[r-(\kappa-1)(\frac{1}{\theta}-\frac{1}{\eta}) ]}j^{-(\frac{1}{\theta}-\frac{1}{\eta})(d-1)}\Big]^p \\
&\lesssim  2^{-L(\frac{1}{\theta}-\frac{1}{\eta} +r)p} L^{-(\frac{1}{\theta}-\frac{1}{\eta})(d-1)p} \\
& \asymp 2^{-(n - (d-1) \log n)(\frac{1}{\theta}-\frac{1}{\eta}+r)p} (n - (d-1) \log n)^{-(\frac{1}{\theta}-\frac{1}{\eta})(d-1)p}\\
& \asymp 2^{-n(\frac{1}{\theta}-\frac{1}{\eta}+r)p} n^{(d-1)(\frac{1}{\theta}-\frac{1}{\eta}+r)p} n^{-(\frac{1}{\theta}-\frac{1}{\eta})(d-1)p} \\
& = 2^{-n(\frac{1}{\theta}-\frac{1}{\eta}+r)p} n^{(d-1)r p}.
\end{split}
\end{align}
The sum in the fourth line converges if \(r > (\kappa-1)(\frac{1}{\theta}-\frac{1}{\eta})\) (see \cite[Lemma 4.2]{MSU26}). Therefore we can choose \(1 < \kappa < \frac{r + \frac{1}{\theta}-\frac{1}{\eta}}{\frac{1}{\theta}-\frac{1}{\eta}} \) and the sum always converges. 

For \(\theta > 1\) we instead set \(m_j = \lceil(j-L)^{-2} 2^L L^{d-1}\rceil\) and apply Lemma \ref{lem:gel_p>1}
\begin{align} \label{eq:gel_calc_1}
\begin{split}
\sum_{j=L+1}^K c_{m_j}^p(\Delta_j\colon \ell_\theta\to\ell_\eta) &\lesssim \sum_{j=L+1}^K c^p_{m_j}(\id\colon \ell_\theta^{C(j)}\to\ell_\eta^{C(j)})2^{-j r p} \\
&\lesssim\sum_{j=L+1}^K \big( (2^j j^{d-1})^{1-1/\theta} (j-L) 2^{-L/2} L^{-(d-1)/2}\big)^{\frac{1/\theta-1/\eta}{1/\theta -1/2}p} 2^{-j r p}\\
& = \big(2^L L^{d-1} \big)^{-\frac{1/\theta-1/\eta}{2/\theta -1}p} \sum_{j=L+1}^K \big((j-L)(2^j j^{d-1})^{1-1/\theta} \big)^{\frac{1/\theta-1/\eta}{1/\theta -1/2}p} 2^{-j r p }\\
&\lesssim  \big(2^{L(\frac{1}{2}-\frac{1}{\theta})} L^{(d-1)(\frac{1}{2}-\frac{1}{\theta})}\big)^{\frac{1/\theta-1/\eta}{1/\theta -1/2}p} 2^{-Lrp} \\
& = 2^{L(\frac{1}{\eta}-\frac{1}{\theta}-r)p} L^{(d-1)(\frac{1}{\eta}-\frac{1}{\theta})p}\\
& \asymp 2^{-(n - (d-1) \log n)(\frac{1}{\theta}-\frac{1}{\eta}+r)p} (n - (d-1) \log n)^{-(\frac{1}{\theta}-\frac{1}{\eta})(d-1)p}\\
& \asymp 2^{-n(\frac{1}{\theta}-\frac{1}{\eta}+r)p} n^{(d-1)(\frac{1}{\theta}-\frac{1}{\eta}+r)p} n^{-(\frac{1}{\theta}-\frac{1}{\eta})(d-1)p} \\
& = 2^{-n(\frac{1}{\theta}-\frac{1}{\eta}+r)p} n^{(d-1)r p}.
\end{split}
\end{align}
This sum converges for \(r > (1- 1/\theta) \big(\frac{1/\theta-1/\eta}{1/\theta -1/2}\big)= \frac{2 (\eta -\theta)(\theta -1)}{\eta (2 -\theta) \theta} \).
Now all that is left is to investigate the setting where \(j \geq K+1\). Here we set \(m_j = 0\), then again by the same arguments as in \eqref{eq:delta_EL}, we have that
\begin{equation}
            c_{0}\big(\sum_{j=K+1}^\infty\Delta_j\big) = E_{\hyp^K}(\ca^{r}_\theta)_{\ca_\eta}.
\end{equation}
In this setting \(E_{\hyp^M}\) scales, by Theorem \ref{thm:E_M_hyp}, like \(2^{-Mr}\).
We can now set
\[K = \Big\lfloor n \Big(1 + \frac{{\frac{1}{\theta}-\frac{1}{\eta}}}{r}\Big) - (d-1) \log n\Big\rfloor,\] and obtain
 \begin{equation}
            c_{0}\big(\sum_{j=K+1}^\infty\Delta_j\big) \asymp  2^{-Kr} \asymp m^{-(r  - \frac{1}{\eta} + \frac{1}{\theta})} \log(m)^{(d-1)r}.
 \end{equation}

 Before concluding the proof of the upper bounds we still need to check that \eqref{eq:decomp_gel_3} actually holds, i.e., that 
\begin{equation} \label{eq:check_mk}
         \sum_{j = 0}^K m_j = \sum_{j = 0}^L C(j) + \sum_{j = L+1}^K m_j \lesssim m.
\end{equation}
The bound for the first sum is given by \eqref{L_upper}. For the second sum we need to differentiate between the two relevant settings. For \(\theta \leq 1\) we have
\begin{equation}
        \sum_{j = L+1}^K \lceil 2^j 2^{(L-j)\kappa}j^{d-1}\rceil \leq K + 2^{L\kappa}\sum_{j = L+1}^K  2^{j(1-\kappa)} j^{d-1}\lesssim m + 2^{L} L^{d-1} \lesssim m
\end{equation} where we used the fact that \( \kappa > 1\) and therefore the sum converges. For \(\theta > 1\) we instead get 
\begin{equation}
        \sum_{j = L+1}^K \lceil(j-L)^{-2} 2^L L^{d-1}\rceil  \leq K + 2^L L^{d-1 }\sum_{j = L+1}^K  (j-L)^{-2} \lesssim m + 2^{L} L^{d-1} \lesssim m.
\end{equation}
    To show lower bounds we will restrict to  \(\Delta_{L+c_d}\). For this layer we have that the number of points on it is at least \(C(L+c_d) \geq 2 m\). For \( \theta \leq 1\) we then have by Lemma \ref{lem:gel_rates_ell},
\begin{align} 
\begin{split}
 c_m(\id:\ca^{r+s}_\theta \to \ca^s_\eta) & \geq c_{m}(\Delta_{L+c_d}\colon \ell_\theta\to\ell_\eta) \\
 & \gtrsim c_{m}(\id\colon \ell_\theta^{C(L+c_d)}\to\ell_\eta^{C(L+c_d)})2^{-r(L+c_d)} \\
&\asymp \bigg(\frac{\log\Big(\frac{2^{L+c_d}(L+c_d)^{d-1}}{m}\Big)}{m}\bigg)^{(\frac{1}{\theta}-\frac{1}{\eta})} 2^{-r(L+c_d)}\\ 
&\asymp \log\Big(\frac{2m}{m}\Big)^{(\frac{1}{\theta}-\frac{1}{\eta})}{m}^{-(\frac{1}{\theta}-\frac{1}{\eta})} 2^{-rL}\\
& \asymp m^{-(r -\frac{1}{\eta} + \frac{1}{\theta})} \log(m)^{(d-1)r}.\\
\end{split}
\end{align}
For \(\theta > 1\) we again use Lemma \ref{lem:gel_p>1} and obtain
\begin{align} 
\begin{split}
c_m(\id:\ca^{r+s}_\theta \to \ca^s_\eta) & \geq c_{m}(\Delta_{L+c_d}\colon \ell_\theta\to\ell_\eta) \\ 
&\gtrsim c_{m}(\id\colon \ell_\theta^{C(L+c_d)}\to\ell_\eta^{C(L+c_d)})2^{-r(L+c_d)} \\
&\asymp \big( (2^{L+c_d} (L+c_d)^{d-1})^{1-1/\theta} m^{-1/2} \big)^{\frac{1/\theta-1/\eta}{1/\theta -1/2}} 2^{-(L+c_d)r}\\
& \asymp  m^{1/\eta-1/\theta} 2^{-Lr}\\
& \asymp m^{-(r -\frac{1}{\eta} + \frac{1}{\theta})} \log(m)^{(d-1)r}.\\
\end{split}
\end{align}
This holds after checking that we correctly applied the minimum in the third line.
Arguing with monotonicity now again concludes the proof.
\end{proof}
The formulation of this result for \(\theta = \eta\) and in particular \(\theta = \eta = 2\) is not proven but follows directly from Theorem \ref{thm:gel_th>eta} above. For \(\theta = \eta = 2\) the condition is \(r>0\). We still included these cases so we do not need to reference Theorem \ref{thm:gel_th>eta} additionally whenever we mean to include the regime where \(\theta = \eta\).

Note that the fact \(\eta \leq 2\) is crucial in the above result since the behavior of the Gelfand numbers changes at \(2\).

\begin{thm}\label{thm:gel_2<th<eta}
     Let \(s \geq 0\), \(r > 0\) and \(2 \leq \theta \leq \eta \leq \infty\), then it holds 
    \begin{equation}
         c_m(\id:\ca^{r+s}_\theta \to \ca^s_\eta) \asymp m^{-r} \log(m)^{(d-1)r}.
    \end{equation}
\end{thm}
\begin{proof}
    We start by noticing that the upper bound is trivial by the decomposition in Equation \eqref{eq:decomp_gel_2} together with the upper bound from Theorem \ref{thm:E_M_hyp}.
    For the lower bound we again consider the hyperbolic layer \(\hyp_{L+c_d}\) which has at least \(2m\) entries.
    Since \(\theta \geq 2\) we get by Lemma \ref{lem:gel_1<p<2} the following
    \begin{align} 
\begin{split}
c_m(\id:\ca^{r+s}_\theta \to \ca^s_\eta) & \geq c_{m}(\Delta_{L+c_d}\colon \ell_2\to\ell_\eta) \\ 
&\gtrsim c_{m}(\id\colon \ell_2^{C(L+c_d)}\to\ell_\eta^{C(L+c_d)})2^{-r(L+c_d)} \\
&\gtrsim \max\big\{(2^{L} L^{d-1})^{1/\eta-1/2}, \sqrt{1/2} \,\big\} 2^{-rL} \\
& \gtrsim m^{-r} \log(m)^{(d-1)r}.
\end{split}
\end{align}

\end{proof}

\begin{thm}\label{thm:gel_th<2<eta}
 For \(s \geq 0\), \(1 \leq \theta \leq 2 \leq \eta \leq \infty\) and \(r > 1 - 1/\theta\) we have
        \begin{equation}
            c_m(\id:\ca_\theta^{r+s} \to \ca_\eta^s) \asymp m^{-(r - \frac{1}{2}  + \frac{1}{\theta})} \log(m)^{(d-1)r}.
        \end{equation}
\end{thm}
\begin{proof}
To obtain the upper bound we can simply apply Theorem \ref{thm:gel_th<eta<2} with \(\eta = 2\), since we have
\begin{equation} \label{eq:turm_ca}
    c_m(\ca_\theta^r)_{\ca_\infty} \leq c_m(\ca_\theta^r)_{\ca_\eta} \leq  c_m(\ca_\theta^r)_{\ca_2},
\end{equation}
for all \(2 \leq \eta \leq \infty\). The restriction on \(r\) then simplifies to \(r > 1 - 1/\theta\).
To show the lower bound we make use of the first inequality in \eqref{eq:turm_ca} and again restrict to \(\hyp_{L+c_d}\) which still has at least \(2m\) entries. We therefore obtain

\begin{align} 
\begin{split}
c_m(\id:\ca^{r+s}_\theta \to \ca^s_\eta) &\gtrsim c_{m}(\id\colon \ell_\theta^{C(L+c_d)}\to\ell_\infty^{C(L+c_d)})2^{-r(L+c_d)} \\
&\gtrsim (2^{L} L^{d-1})^{1-1/\theta} m^{-1/2} 2^{-rL} \sqrt{1/2} \\
& \gtrsim m^{-(r - 1/2 + 1/\theta)} \log(m)^{(d-1)r},
\end{split}
\end{align}
where we used the fact that the second term dominates the maximum from Lemma \ref{lem:gel_1<p<2}, since
\[C(L+c_d)^{-1/\theta} \asymp m^{-1/\theta} \leq m^{1/2-1/ \theta} \lesssim C(L+c_d)^{1-1/\theta} m^{-1/2}.\]
  For \(\theta = 1\) we need to use the modified result to get the lower bound.
\end{proof}

\section{Gelfand numbers in \(L_q\)} \label{sec:gel_Lq}

In this section we develop some results for Gelfand numbers of weighted Wiener spaces in Lebesgue spaces. 
We always have the trivial bounds by Hausdorff-Young together with the \(L_2\)-bound where the upper and lower bounds are flipped between the regimes where \(q>2\) and \(q<2\).

For \(q \leq 2\) we are already in a good situation because the lower bound by Hausdorff-Young will now be the same as the upper bound by \(L_2\) due to the fact that the Gelfand numbers in Theorem \ref{thm:gel_th<2<eta} are independent of \(\eta \geq 2 \).

\begin{thm} \label{thm:gel_L_q<2}
    Let \(1 \leq \theta \leq 2\) and \(1 \leq q \leq 2\), then for \(r > 1 - 1/\theta\) we have
    \begin{equation}
        c_m(\ca_\theta^r)_{L_q} \asymp m^{-(r -\frac{1}{2} + \frac{1}{\theta})} \log(m)^{(d-1) r}.
    \end{equation}
   
\end{thm}
\begin{proof}
    We follow the strategy outlined above and use the rates for Gelfand numbers between Wiener spaces from Theorem \ref{thm:gel_th<eta<2}.
    We begin with the upper bound which is given by the \(L_2\) embedding.
    \begin{equation}
        c_m(\ca_\theta^r)_{L_q} \leq c_m(\ca_\theta^r)_{\ca_{2}} \asymp m^{-(r -  \frac{1}{2} + \frac{1}{\theta})}\log(m)^{(d-1) r}.
    \end{equation}
    The corresponding lower bound is now given by Hausdorff-Young, Proposition \ref{prop:Haus_young} together with Theorem \ref{thm:gel_th<2<eta}. Note that this is the correct bound to apply since we restricted \(\theta \leq 2 \leq q'\). We obtain 
        \begin{equation}
        c_m(\ca_\theta^r)_{L_q} \geq c_m(\ca_\theta^r)_{\ca_{q'}} \asymp m^{-(r - \frac{1}{2} + \frac{1}{\theta})}\log(m)^{(d-1) r}.
    \end{equation}
\end{proof}

 For \(\theta > 2\) the upper bounds still hold but we cannot show lower bounds with the same rate. Instead we get 
 \begin{equation}
     c_m(\ca_\theta^r)_{L_q} \gtrsim m^{-r}\log(m)^{(d-1) r},
 \end{equation}  for \( 2 \leq \theta \leq q'\).
 For \( q' \leq \theta\) we instead have the trivial bound
  \begin{equation}
     c_m(\ca_\theta^r)_{L_q} \gtrsim m^{-(r -1 + \frac{1}{q} + \frac{1}{\theta})}\log(m)^{(d-1) r}.
 \end{equation} 

Let us now transition to the setting where \(q \geq 2\). Here a simple application of the Hausdorff-Young inequality is not sufficient, instead we need an additional lemma.

\begin{lem}[{\cite[Ch. 2 Lemma 2.1]{Tem93_book}}] \label{lem:tem93}
    Let \(1 < p < q< \infty\) and \(f \in L_p\), then it holds
    \begin{equation}
        \Vert f \Vert_{L_q} \lesssim \Big( \sum_{\bj \in \Z} \big\Vert \sum_{\bk \in I_\bj} \exp(2 \pi \mathrm{i} \bk \bx) \hat{f}(\bk) \big\Vert_{L_p}^q 2^{\Vert \bj \Vert_{1}(q/p-1)} \Big)^{1/q}.
    \end{equation}
\end{lem}

\begin{thm}\label{thm:gel_L_q>2}
    Let \(0 <\theta \leq 2 \leq q < \infty\) then it holds for \(r > \frac{1}{2} - \frac{1}{q} + \big(1 - \frac{1}{\theta}\big)_+\)
    \begin{equation}
      c_m(\ca_\theta^r)_{L_q} \asymp m^{-(r - 1 + \frac{1}{q} + \frac{1}{\theta})}\log(m)^{(d-1)(r - \frac{1}{2} +\frac{1}{q})}.
    \end{equation}
    \end{thm}
\begin{proof}
    
    Let \(\theta \leq 1\). By the proof of Theorem \ref{thm:gel_th<eta<2} we know that there exists a closed linear subspace \(M \subset \ca^r_\theta\) with \(\operatorname{codim} M \lesssim m\) and \(\hyp^L \cap M = \{0\}\) for \(L\) as in \eqref{L} that has codimension \(m_l\) in \(\hyp_l\) for \(l = L+1, \ldots, K\), where \(m_l = \lceil 2^l 2^{(L-l)\kappa}l^{d-1} \rceil\) such that \(M\) recovers the bound for Gelfand numbers between Wiener spaces in Theorem \ref{thm:gel_th<eta<2}. We now fix such a subspace \(M\) and apply Lemma \ref{lem:tem93} to obtain for \(f \in M \cap \hyp^K\), that 
 \begin{align} \label{eq:gel_Lq_improv_rates}
    \begin{split}
          \sup_{\Vert f\Vert_{\ca_\theta^r}\leq 1} \Vert f \Vert_{L_q} & \lesssim \sup_{\Vert f\Vert_{\ca_\theta^r}\leq 1} \Big( \sum_{\Vert \bj \Vert_1 \in \n} \Big\Vert \sum_{\bk \in I_\bj} \exp(2 \pi \mathrm{i} \bk \bx) \hat{f}(\bk) \Big\Vert_{L_2}^q 2^{\Vert \bj \Vert_{1}(q/2-1)} \Big)^{1/q} \\
         & = \sup_{\Vert f\Vert_{\ca_\theta^r}\leq 1} \Big( \sum_{l = L+1}^K 2^{lq(1/2-1/q)} \sum_{\Vert \bj \Vert_1 = l} \big\Vert f_{I_\bj} \big\Vert_{\ca_{2}}^q  \Big)^{1/q} \\
         & \leq \sup_{\Vert f\Vert_{\ca_\theta^r}\leq 1} \Big( \sum_{l = L+1}^K 2^{lq(1/2-1/q)} \big\Vert f_{\hyp_l} \big\Vert_{\ca_{2}}^q  \Big)^{1/q} \\
    \end{split}
\end{align}
Then we can apply Lemma \ref{lem:gel_rates_ell} in the same way as in Theorem \ref{thm:gel_th<eta<2} and get the same blockwise estimates. For a subspace \(M \subset \ca^r_\theta\) as above and \(f \in M \cap \hyp^K\) it holds
\begin{equation} \label{eq:gel_block_rates_theta<1}
    \sup_{\Vert f\Vert_{\ca_\theta^r}\leq 1} \big\Vert f_{\hyp_l} \big\Vert_{\ca_{2}} \lesssim (l - L)\kappa 2^{-l[r-(\kappa-1)(\frac{1}{\theta}-\frac{1}{2}) ]} 2^{-L \kappa(\frac{1}{\theta}-\frac{1}{2})} l^{-(\frac{1}{\theta} - \frac{1}{2})(d-1)}.
\end{equation}
Since our subspace \(M \subset \ca^r_\theta\) has \(\operatorname{codim} M \lesssim m\) and  therefore (up to some constant) is admissible for \(c_m\) we can now combine \eqref{eq:gel_Lq_improv_rates} with \eqref{eq:gel_block_rates_theta<1}. We get
\begin{align}
    \begin{split}
        c_m(\ca^r_\theta)_{L_q} & \lesssim  \Big( \sum_{l = L+1}^K 2^{lq(\frac{1}{2}-\frac{1}{q})} (l - L)^q 2^{-lq[r-(\kappa-1)(\frac{1}{\theta}-\frac{1}{2}) ]} 2^{-L \kappa q(\frac{1}{\theta}-\frac{1}{2})} l^{-q(\frac{1}{\theta} - \frac{1}{2})(d-1)}  \Big)^{1/q} \\
        & \leq 2^{-L \kappa (\frac{1}{\theta}-\frac{1}{2})} \Big( \sum_{l = L+1}^K  (l - L)^q 2^{-lq[r-\frac{1}{2}+\frac{1}{q}-(\kappa-1)(\frac{1}{\theta}-\frac{1}{2}) ]} l^{-q(\frac{1}{\theta} - \frac{1}{2})(d-1)}  \Big)^{1/q}\\
        & \lesssim 2^{-L \kappa (\frac{1}{\theta}-\frac{1}{2})} 2^{-L[r- \frac{1}{2}+\frac{1}{q}-(\kappa-1)(\frac{1}{\theta}-\frac{1}{2}) ]} L^{-(\frac{1}{\theta} - \frac{1}{2})(d-1)} \\
        & \leq 2^{-L(r- 1+ \frac{1}{q}+\frac{1}{\theta})} L^{-(\frac{1}{\theta} - \frac{1}{2})(d-1)}\\
        & \lesssim m^{-(r- 1+ \frac{1}{q}+\frac{1}{\theta})} \log(m)^{(d-1)(r -1+ \frac{1}{q} + \frac{1}{\theta} - \frac{1}{\theta} + \frac{1}{2})}\\
        & = m^{-(r- 1+ \frac{1}{q}+\frac{1}{\theta})} \log(m)^{(d-1)(r -\frac{1}{2}+ \frac{1}{q} )},
    \end{split}
\end{align}
where we needed \(r- \frac{1}{2}+\frac{1}{q} >(\kappa-1)(\frac{1}{\theta}-\frac{1}{2}) \) for the sum to converge by the same argument as in the proof of Theorem \ref{thm:gel_th<eta<2}.
Additionally this requires us to bound \(E_{\hyp^K}(\ca^r_\theta)_{L_q} \) with at worst the  same error bound so we can drop the additional restriction \(f \in \hyp^K\). This follows from Hausdorff-Young, Proposition \ref{prop:Haus_young} together with the upper bound from Theorem \ref{thm:E_M_hyp}. Since  \(\theta < q'\) here we set \(K = \frac{1}{r}\Big\lfloor\big(r -1 +\frac{1}{q} + \frac{1}{\theta}\big) n -  \big(r - \frac{1}{2} + \frac{1}{q} \big) (d-1) \log n\Big\rfloor\) and get
\begin{equation}
    E_{\hyp^K}(\ca^r_\theta)_{L_q} \leq E_{\hyp^K}(\ca^r_\theta)_{\ca_{q'}} \lesssim 2^{-n(r - 1 +\frac{1}{q} + \frac{1}{\theta})} n^{(d-1)(r- \frac{1}{2} +\frac{1}{q})} \leq m^{-(r - 1 +\frac{1}{q} + \frac{1}{\theta}) } \log(m) ^{(d-1)(r- \frac{1}{2} +\frac{1}{q})}.
\end{equation}

Let us now turn to the regime where \(1 < \theta \leq 2\). We again use \eqref{eq:gel_Lq_improv_rates} but instead of \eqref{eq:gel_block_rates_theta<1} we need the block result given by Lemma \ref{lem:gel_p>1}. We restrict \(M\) to be of the form given by the choice of \(m_l = \lceil(l-L)^{-2} 2^L L^{d-1}\rceil\) there. From the proof of Theorem \ref{thm:gel_th<eta<2} we then get (for \(\eta = 2\)) that
\begin{equation}  \label{eq:gel_block_rates_theta>1}
    \sup_{\Vert f\Vert_{\ca_\theta^r}\leq 1} \Vert f_{\hyp_l} \Vert_{\ca_{2}} \lesssim (l-L) (2^l l^{d-1})^{1-1/\theta} 2^{-rl}  2^{-L/2} L^{-(d-1)/2} .
\end{equation}

For \(\theta \leq q'\) we again have the bound on \(E_{\hyp^K}(\ca^r_\theta)_{L_q}\) as before, but for \(\theta \geq q'\) it is already sufficient to set \(K =n\), then we have by Theorem \ref{thm:E_M_hyp} that
\begin{equation}
    E_{\hyp^K}(\ca^r_\theta)_{L_q} \leq E_{\hyp^K}(\ca^r_\theta)_{\ca_{q'}} \lesssim 2^{-n(r - 1 +\frac{1}{q} + \frac{1}{\theta})} n^{(d-1)(1- \frac{1}{q}- \frac{1}{\theta})} \leq m^{-(r - 1 +\frac{1}{q} + \frac{1}{\theta}) } \log(m) ^{(d-1)(r- \frac{1}{2})}.
\end{equation}
Combining this with \eqref{eq:gel_Lq_improv_rates} and \eqref{eq:gel_block_rates_theta>1} yields
\begin{align}
    \begin{split}
         c_m(\ca^r_\theta)_{L_q} & \lesssim  \Big( \sum_{l = L+1}^K 2^{lq(\frac{1}{2}-\frac{1}{q})} (l - L)^q (2^l l^{d-1})^{q(1-\frac{1}{\theta})} 2^{-rlq}  2^{-qL/2} L^{-q(d-1)/2}  \Big)^{1/q} \\
         &\leq 2^{-L/2} L^{-(d-1)/2}  \Big( \sum_{l = L+1}^K (l - L)^q  2^{-lq(r - \frac{3}{2}+\frac{1}{q}+\frac{1}{\theta})} l^{q(d-1)(1-\frac{1}{\theta})} \Big)^{1/q} \\
         & \lesssim  2^{-L(r - 1+\frac{1}{q}+\frac{1}{\theta})} L^{(d-1)(\frac{1}{2}-\frac{1}{\theta})} \\
         & \lesssim m^{-(r - 1+\frac{1}{q}+\frac{1}{\theta})} \log(m)^{(d-1)(r - \frac{1}{2} +\frac{1}{q})}.
    \end{split}
\end{align}

   To get the lower bounds we will combine multiple results. First by Definition \ref{def:sfunc} (S2) (and Proposition \ref{prop:all_s_num}) together with Theorem \ref{thm:gel_th>eta} we have for any \(s > r+ (\frac{1}{\theta}-\frac{1}{2}) \)
   \begin{equation}
      c_m(\ca^r_\theta)_{L_q} \geq\frac{ c_{2m}(\ca^s_2)_{L_q}}{c_m(\ca_2^s)_{\ca^r_\theta} } \asymp m^{s-r + \frac{1}{2} - \frac{1}{\theta}} c_{m}(\ca^s_2)_{L_q} \log(m)^{(r-s)(d-1)}.
   \end{equation}
We now want to use the embedding \(B^{s}_{2,2} \hookrightarrow \ca^s_2\), see, e.g., \cite[Lemma 3.12]{DTU18}. Now we can apply known rates for Besov spaces with dominating mixed smoothness from \cite[Theorem 2.2]{Ngu17}
\begin{align}
\begin{split}
c_m(\ca_\theta^r)_{L_q} & \gtrsim c_{m}(\ca^s_2)_{L_q}\, m^{s-r + \frac{1}{2} - \frac{1}{\theta}}  \log(m)^{(r-s)(d-1)} \\
 &\geq c_{m}(B^s_{2,2})_{L_q}\, m^{s-r +  \frac{1}{2} - \frac{1}{\theta}} \log(m)^{(r-s)(d-1)} \\
& \asymp  m^{-(s-(\frac{1}{2} - \frac{1}{q}))} \log(m)^{(d-1)(s-(\frac{1}{2} - \frac{1}{q}))} m^{s-r +  \frac{1}{2} - \frac{1}{\theta}} \log(m)^{(r-s)(d-1)} \\
& = m^{-(r -1 + \frac{1}{q} + \frac{1}{\theta})} \log(m)^{(d-1)(r -  \frac{1}{2} + \frac{1}{q})}.  \\
    \end{split}
\end{align}
\end{proof}

Note that the condition on \(r\) is much more restrictive than what the embedding \(\ca^r_\theta \hookrightarrow L_q\) requires. The lower bound shown here is much sharper than the one obtained by embedding into the one-dimensional Besov spaces as is done in \cite[Remark 1]{MPU24}. The lower bounds are consistent with ones that could be obtained by using results on Kolmogorov numbers of Sobolev spaces with bounded mixed derivatives from, e.g., \cite[Theorem 4.13]{DTU18} together with the Gelfand-Kolmogorov duality and a lifting result, see, e.g., \cite[Section 8]{ByUl17}. 
For \(\theta > 2\) the lower bounds fails since \(c_m(\ca_2^s)_{\ca^r_\theta}\) changes its behavior at \(\theta = 2\).

For \(q = \infty\) we have a different logarithmic rate. 
\begin{cor} \label{cor:gel_L_infty}
    Let \(0 < \theta \leq 2\) and \(r > \big(1 - \frac{1}{ \theta}\big)_+\), then it holds
    \begin{equation}
        c_m(\ca^r_\theta)_{L_\infty} \asymp m^{-(r - 1 + \frac{1}{\theta})}\log(m)^{(d-1)r}.
    \end{equation}
\end{cor}
\begin{proof}
    The upper bound is a simple application of the Hausdorff-Young inequality Proposition \ref{prop:Haus_young} together with the upper bounds from Theorem \ref{thm:gel_th<eta<2} and Theorem \ref{thm:gel_th>eta}, i.e.,
    \begin{equation}
        c_m(\ca^r_\theta)_{L_\infty} \leq c_m(\ca^r_\theta)_{\ca_1} \leq m^{-(r - 1 + \frac{1}{\theta})}\log(m)^{(d-1)r}.
    \end{equation}
    For the lower bound we more closely follow the argument from Theorem \ref{thm:gel_L_q>2}. We again use the embedding \(B^s_{2,2} \hookrightarrow \ca^s_2\). We have by \cite[Theorem 2.8]{Ngu17} for \( s > r+ \frac{1}{\theta} > 1 \) that it holds
    \begin{align}
        \begin{split}
       c_m(\ca_\theta^r)_{L_\infty} & \gtrsim  c_m(\ca_2^s)_{L_\infty} \, m^{s-r + \frac{1}{2} - \frac{1}{\theta}}  \log(m)^{(d-1)(r-s)} \\
    & \geq c_m(B^s_{2,2})_{L_\infty} \, m^{s-r + \frac{1}{2} - \frac{1}{\theta}}  \log(m)^{(d-1)(r-s)} \\
       & \asymp m^{-(s-\frac{1}{2})} \log(m)^{(d-1)s}  m^{s-r + \frac{1}{2} - \frac{1}{\theta}}  \log(m)^{(d-1)(r-s)}  \\
       & = m^{-( r - 1 + \frac{1}{\theta})}  \log(m)^{(d-1)r}. \\
        \end{split}
    \end{align}
  \end{proof}

\section{Consequences for sampling recovery in \(L_q\)} \label{sec:samp_Lq}
In this section we want to combine the results from \cite{MSU26} that give upper bounds on the sampling numbers with the lower bounds for Gelfand numbers developed in this paper.
We will briefly restate the upper bounds and begin with the central inequality that puts the sampling numbers and the best \(m\)-term approximation in relation. As mentioned in the introduction, the first result of this type was shown in \cite{JUV23} and then later generalized in \cite{ABG26, Kri23, MPU24, MSU26}.

\begin{prop}[{\cite[Theorem 3.7]{MSU26}}]\label{prop:samp_A1}
Let \(X \hookrightarrow C(\tor^d)\) be a function space, then for any \(M, m \in \N\) and \(2 \leq q \leq \infty\) it holds
\begin{equation} 
 \varrho_{\lceil m d \log(m)^2 \log(M)\rceil}(X)_{L_q} \lesssim m^{1/2-1/q}\Big(m^{-1/2} \sigma_m\big(X\big)_{  {\ca} } + E_{[-M,M]^d}\big(X\big)_{L_\infty}\Big)\, .
\end{equation}

\end{prop}

Also from that paper we have a sharp bound of the best \(m\)-term approximation between Wiener spaces.

\begin{prop}[{\cite[Theorem 4.5]{MSU26}}]\label{prop:bestm_ca}
For \(r > \big(\frac{1}{\eta} - \frac{1}{\theta}\big)_+\) and   {\( \theta, \eta \in (0,\infty]\)} it holds
\begin{equation}
\sigma_m(\ca^r_\theta)_{\ca_\eta} \asymp 
m^{\frac{1}{\eta} - \frac{1}{\theta}-r}\log(m)^{(d-1)r}.
\end{equation}
\end{prop}

In addition we also need an upper bound on the error of best trigonometric approximation in \(L_\infty\). The upper bound on that quantity follows directly from Theorem \ref{thm:E_M} together with the embedding \(\ca \hookrightarrow L_\infty\). The lower bound can be found in \cite[Proposition 4.6]{ABG26}.

\begin{prop}\label{prop:E_M/L_q}
    For \(0 <\theta \leq \infty\) and \( r > (1 - \frac{1}{\theta})_+\) it holds
    \begin{equation}
         E_{[-M,M]^d}(\ca^{r}_\theta)_{L_\infty} \asymp 
        \begin{cases}
            M^{-r} & \theta \leq 1 \\
            M^{-(r -1 + \frac{1}{\theta})} & \theta > 1.
        \end{cases}
    \end{equation}
\end{prop}

We can now assemble the upper bounds on the sampling numbers from the previous results together with the lower bounds via Gelfand numbers from Section \ref{sec:gel_Lq}.
\begin{thm}\label{thm:samp_main}
    For \(2 \leq q < \infty\), \(0 <\theta \leq 2\) and \(r > \frac{1}{2} - \frac{1}{q} +(1 - \frac{1}{\theta})_+\) we get
    \begin{equation}
          m^{-(r - 1 + \frac{1}{q} + \frac{1}{\theta})}\log(m)^{(d-1)(r -  \frac{1}{2} + \frac{1}{q})}\lesssim\varrho_m(\ca^r_\theta)_{L_q} \lesssim m^{-(r - 1 +  \frac{1}{q} + \frac{1}{\theta})} \log(m)^{(d-1)r +3 (r - 1  + \frac{1}{q}+ \frac{1}{\theta}) }.
    \end{equation}
\end{thm}
\begin{proof}
    The lower bound is simply Theorem \ref{thm:gel_L_q>2} together with Lemma \ref{lem:gel_samp}. 
    To show the upper bound we employ Proposition \ref{prop:samp_A1}. From there we have
\begin{equation} \label{eq:app_samp_bestm}
     \varrho_{\lceil m d \log(m)^2 \log(M)\rceil}(\ca^r_\theta)_{L_q} \lesssim m^{-1/q}\sigma_m(\ca^r_\theta)_{\ca} +  m^{1/2-1/q} E_{[-M,M]^d}(\ca^r_\theta)_{L_\infty}.
\end{equation}
    The first term of the right hand side of this inequality can be dealt with by Proposition \ref{prop:bestm_ca}. This yields
    \begin{equation}
        \sigma_m(\ca^r_\theta)_{\ca} \asymp m^{-(r - 1 + \frac{1}{\theta})}\log(m)^{(d-1)r}.
    \end{equation}
    
    We can now employ Proposition \ref{prop:E_M/L_q} to deal with the second term in Equation \eqref{eq:app_samp_bestm}. For \(\theta \geq 1 \) we set \(M = \Big\lceil m^{\frac{r + 1/\theta - 1/2}{r -1 +1/\theta}} \Big\rceil\) and for \(\theta \leq 1 \) as \(M = \Big\lceil m^{\frac{r  + 1/\theta - 1/2}{r }} \Big\rceil\) instead. This then yields 
    \begin{equation}
        E_{[-M,M]^d}(\ca^r_\theta)_{L_\infty} \asymp m^{-(r  + 1/\theta - 1/2)}.
    \end{equation}
    
    All that is left to do is to rephrase the left hand side in terms of \(m\) instead of \(\lceil m d \log(m)^2 \log(M)\rceil\) which yields the assertion.
\end{proof}
Note that for \(\theta > 2\) the upper bounds will still hold as the results for the best \(m\)-term approximation have no comparable restrictions on \(\theta\).
For \(q = \infty\) we get the slightly larger lower bound from Corollary \ref{cor:gel_L_infty} instead.
\begin{cor}\label{cor:samp_L_infty}
        For \(0 <\theta \leq 2\) and \(r > (1 - \frac{1}{\theta})_+\) we get
    \begin{equation}
          m^{-(r - 1  + \frac{1}{\theta})}\log(m)^{(d-1)r}\lesssim\varrho_m(\ca^r_\theta)_{L_\infty} \lesssim m^{-(r - 1 + \frac{1}{\theta})} \log(m)^{(d-1)r +3 (r - 1 + \frac{1}{\theta}) }.
    \end{equation}
\end{cor}
\begin{proof}
    The lower bound follows from Corollary \ref{cor:gel_L_infty} together with Lemma \ref{lem:gel_samp}. The upper bound follows from the same machinery as before with \(q = \infty\) in Proposition \ref{prop:samp_A1}.
\end{proof}

\begin{rem}
    For the precise situation where \(q = \infty\) and \(\theta = 1\), Smolyak Sampling type algorithms give a better upper bound on the non-linear sampling numbers, see \cite[Theorem 5.1 and Remark 6.4]{KLT23}. Here we have for all \(r> 0\) that 
    \begin{equation}
        \varrho_m(\ca^r_1)_{L_\infty} \asymp m^{-r } \log(m)^{(d-1)r }.
    \end{equation}
    This comparison was first made in \cite[Remark 22]{KPUU23}, where it was observed that Smolyak Sampling gives the best known rates for the linear sampling numbers for \(2 < q \leq \infty\), beating rates given by algorithms with random points. 
    For \(\theta > 1\) the results from \cite{KLT23} give worse dimension dependent logarithmic rates and for \(q' > \theta\) the main rates between linear and non-linear sampling numbers diverge.
\end{rem}

\begin{rem}
    The gap between the lower and upper bounds in Theorem  \ref{thm:samp_main} for \(2 < q < \infty\) is significant and not something that could be easily fixed in the framework used in this paper. From \cite[Theorem 4.3]{MNPSU26} we know that the bound in Proposition \ref{prop:samp_A1} is almost sharp and that at least the logarithmic oversampling in \(M\) is needed. 
    The results for the best \(m\)-term approximation from \cite{MSU26} as well as the results for Gelfand numbers presented in this paper are sharp. However, even if we discard the non-dimension-dependent rates and compare \(c_m(\ca^r_\theta)_{L_q}\) with \(m^{- \frac{1}{q}} \sigma_m(\ca^r_\theta)_{\ca}\) we can see that the Gelfand numbers decay faster by a factor of \(\log(m)^{(d-1)(\frac{1}{2} - \frac{1}{q})}\) compared to the best \(m\)-term approximation. 
    Therefore either the Gelfand numbers as a lower bound on the sampling numbers or the construction of the upper bound admits a genuine dimension-dependent logarithmic gap when bounding the sampling numbers of weighted Wiener spaces in \(L_q\). 
    Note that this gap does not appear when measuring the sampling numbers in \(L_\infty\).
\end{rem}

For \(q \leq 2\) we instead again get that the rates for \(L_2\) are the correct ones.
\begin{cor}
 Let \(1 \leq q \leq 2\), \(1 \leq \theta \leq 2\) and \(r > 1 - \frac{1}{\theta}\) then we have
 \begin{equation}
     m^{-(r -\frac{1}{2} + \frac{1}{\theta})} \log(m)^{(d-1) r} \lesssim \varrho_m(\ca^r_\theta)_{L_q} \lesssim m^{-(r - \frac{1}{2} + \frac{1}{\theta})} \log(m)^{(d-1)r +3 (r - \frac{1}{2} + \frac{1}{\theta}) }.
 \end{equation}
\end{cor}
\begin{proof}
The lower bound is again trivial by Theorem \ref{thm:gel_L_q<2} together with Lemma \ref{lem:gel_samp}. The upper bound is Theorem \ref{thm:samp_main} for \(L_2\) together with the embedding \(L_2 \hookrightarrow L_q\) for \(q \leq 2\).
\end{proof}

\nocite{Kas77E,Glu84E}

\paragraph{Acknowledgments}
   The author is supported by the ESF, being co-financed by the European Union and from tax revenues on the basis of the budget adopted by the Saxon State Parliament.

\printbibliography

\end{document}